\documentclass{ifacconf}

\usepackage{graphicx}      % include this line if your document contains figures
\usepackage{natbib}        % required for bibliography
\usepackage{bm}
\usepackage{amsmath}
\usepackage{ifplatform}
\usepackage{amssymb}
\usepackage{booktabs}
\usepackage{cancel}
\usepackage{marvosym}
\usepackage[mode=buildnew]{standalone}% requires -shell-escape
\usepackage{tikz}
\usepackage{pgfplots}
\usepackage{algorithm}
\usepackage{xspace}
\usepackage{siunitx}
\usepackage[noend]{algpseudocode}
\algnewcommand\algorithmicforeach{\textbf{for each}}
\algdef{S}[FOR]{ForEach}[1]{\algorithmicforeach\ #1\ \algorithmicdo}
\usetikzlibrary{positioning, fit, arrows.meta, shapes}
\newcommand{\R}{\ensuremath\mathbb{R}}
\newcommand{\state}{\bm{z}}
\newcommand{\stateApprox}{\tilde{\state}}
\newcommand{\stateSpace}{\mathcal{Z}}

\newcommand{\redStateSpace}{\mathcal{Z}_r}
\newcommand{\redState}{\bar{\state}}
\newcommand{\redStateApprox}{\tilde{\redState}}
\newcommand{\deltaRedState}{\Delta\redState}
\newcommand{\deltaRedStateApprox}{\Delta\redStateApprox}
\newcommand{\redStateDim}{r}
\newcommand{\stateDim}{N}

\newcommand{\param}{\bm{\mu}}
\newcommand{\paramSpace}{\mathcal{P}}
\newcommand{\paramFuncSpace}{\mathcal{M}}
\newcommand{\paramSet}{\mathcal{M}_\nSnaps}
\newcommand{\paramDim}{\ell}

\newcommand{\timeInt}{\mathcal{T}}
\newcommand{\nTime}{\eta}
\newcommand{\indSpace}{\mathcal{I}}
\newcommand{\snaps}{\bm{Z}}
\newcommand{\nSnaps}{\kappa}
\newcommand{\nRed}{n_\text{red}}

\newcommand{\redMatrix}{\bm{V}}
\newcommand{\redBasis}{\bm{v}}

\newcommand{\Flow}{\bm{F}}
\newcommand{\transition}{\bm{\Phi}}
\newcommand{\difference}{\bm{\phi}}
\newcommand{\reduction}{\bm{V}^T}
\newcommand{\reconstruction}{\bm{V}}

\newcommand{\ODE}{\textsf{ODE}\xspace}
\newcommand{\ODEs}{\textsf{ODEs}\xspace}
\newcommand{\FE}{\textsf{FE}\xspace}
\newcommand{\MBS}{\textsf{MBS}\xspace}
\newcommand{\PDE}{\textsf{PDE}\xspace}
\newcommand{\PDEs}{\textsf{PDEs}\xspace}
\newcommand{\LSTM}{\textsf{LSTM}\xspace}
\newcommand{\LSTMs}{\textsf{LSTMs}\xspace}

\newcommand{\RNNs}{\textsf{RNNs}\xspace}
\newcommand{\MOR}{\textsf{MOR}\xspace}
\newcommand{\POD}{\textsf{POD}\xspace}
\newcommand{\SVD}{\textsf{SVD}\xspace}
\newcommand{\ML}{\textsf{ML}\xspace}
\newcommand{\weights}{\bm{w}}
\newcommand{\Weights}{\bm{W}}
\newcommand{\bias}{\bm{b}}
\newcommand{\NNinput}{\bm{x}}
\newcommand{\NNoutput}{\bm{y}}
\newcommand{\NNhiddenState}{\bm{h}}
\newcommand{\NNinputDim}{n_{x}}
\newcommand{\NNhiddenDim}{n_h}
\newcommand{\NNoutputDim}{n_{y}}
\newcommand{\timeHorizon}{n_w}
\newcommand{\loss}{\mathcal{L}}
\newcommand{\dataSet}{\bm{D}}
\newcommand{\sigmoid}{\bm{\sigma}}

\newcommand{\eRec}{s_{\text{rec}}}
\newcommand{\eRegr}{s_{\text{regr}}}
\newcommand{\eApprox}{s_{\text{approx}}}
\newcommand{\eDist}{e_{\text{dist}}}
\newcommand{\eDistMax}{e_{\text{dist}}^{\text{max}}}
\newcommand{\tRatio}{\Delta t_\text{r}}

\newcommand{\nodeDisp}{\bm{q}}

\newcommand{\nodeDim}{25411}
\newcommand{\dofs}{76233}

\begin{document}
\begin{frontmatter}

\title{Real-time Human Response Prediction Using a Non-intrusive Data-driven Model Reduction Scheme\thanksref{footnoteinfo}\thanksref{footnoteinfo2}
} 

\thanks[footnoteinfo]{Funded by Deutsche Forschungsgemeinschaft (DFG, German Research Foundation) under Germany's Excellence Strategy - EXC 2075 – 390740016. We acknowledge the support by the Stuttgart Center for Simulation Science (SimTech).}

\thanks[footnoteinfo2]{This work has been submitted to IFAC for possible publication.}

\author[First]{J. Kneifl} 
\author[First,Second]{J. Hay} 
\author[First]{J. Fehr}

\address[First]{Institute of Engineering and Computational Mechanics, \\University of Stuttgart, Pfaffenwaldring 9, 70569 Stuttgart, Germany. \\(e-mail: \{jonas.kneifl, joerg.fehr\}@itm.uni-stuttgart.de)}
%\address[Second]{Institute of Engineering and Computational Mechanics, University of Stuttgart, Stuttgart, Germany. (e-mail: julian.hay@itm.uni-stuttgart.de)}
\address[Second]{ZF Friedrichshafen AG,	Safe Mobility Simulation, \\88046 Friedrichshafen, Germany.
	(e-mail: julian.hay@zf.com)}

\begin{abstract}                % Abstract of not more than 250 words.
Recent research in non-intrusive data-driven model order reduction (\MOR) enabled accurate and efficient approximation of parameterized ordinary differential equations (\ODEs). However, previous studies have focused on constant parameters, whereas time-dependent parameters have been neglected. 
The purpose of this paper is to introduce a novel two-step \MOR scheme to tackle this issue. In a first step, classic \MOR approaches are applied to calculate a low-dimensional representation of high-dimensional \ODE solutions, i.e. to extract the most important features of simulation data. Based on this representation, a long short-term memory (\LSTM) is trained to predict the reduced dynamics iteratively in a second step. This enables the parameters to be taken into account during the respective time step. The potential of this approach is demonstrated on an occupant model within a car driving scenario. The reduced model's response to time-varying accelerations matches the reference data with high accuracy for a limited amount of time. Furthermore, real-time capability is achieved.
Accordingly, it is concluded that the presented method is well suited to approximate parameterized \ODEs and can handle time-dependent parameters in contrast to common methods.
\end{abstract}

\begin{keyword}
Model Reduction, Machine Learning, Occupant Safety, Human Body Modeling, Parameterized Ordinary Differential Equations, Long Short-Term Memory.
\end{keyword}

\end{frontmatter}
%===============================================================================

\section{Introduction}
The development of future vehicle safety systems requires knowledge of the driver's movement in the pre-crash phase to engage the safety system optimally and to reduce the resulting injuries. 
For multi-parameter simulations and on-board applications real-time predictions of the human response are necessary. 
However, with common high-fidelity human body models, this is not feasible because they rely on finite element (\FE) or elastic multi- body system (\MBS) formulations.

As a result, a high demand for efficient but accurate surrogate models arises. In the case of nonlinear black-box models, data-based non-intrusive \emph{model order reduction} (\MOR) schemes have proven their potential. Furthermore, they can benefit from the development, research and validation effort included in the complex \MBS and \FE systems.
One possible approach is to reduce the dimension of the system first using \emph{proper orthogonal decomposition} (\POD) and then learn the reduced dynamics using regression algorithms from the field of \emph{machine learning} (\ML). 
For example, \cite{HesthavenUbbiali18} used neural networks to predict the reduced dynamics of parametrized partial differential equations (\PDEs).

Most of the literature in this field focuses on predictions for fixed parameters like physical or geometrical properties of the system, e.g. material properties or initial and boundary values. However, time-dependent parameter trajectories (vehicle accelerations) must be taken into account for the occupant model investigated in this paper. Hence, current ML methods reach their limits regarding approximation quality. One approach to tackle time-dependent parameter dependencies is suggested by~\cite{Zhuang_2021}, who apply model order reduction based on a Runge-Kutta neural network. 
In contrast, a long short-term memory network (\LSTM) is implemented in this work enhancing the iterative approach presented by the authors in \cite{KneiflGrunertFehr21} to capture the time-dependencies. Such networks use recurrent loops to discover temporally encoded information and are commonly used for sequential tasks like language translation or speech recognition. Moreover, they have been used for data-driven time series forecasting in a reduced space in \cite{VlachasEtAl2018}.

In summary, two main problems are addressed in this paper. First, how occupant behavior can be efficiently predicted in the pre-crash phase, and second, how parameter time dependencies can be incorporated into a surrogate model.
The paper itself is structured as follows: section~\ref{sec:method} introduces the occupant model, gives an overview of the required theoretical knowledge and explains the procedure to conduct the non-intrusive \MOR approach combining \POD and \LSTMs. In section~\ref{sec:implementation} the implementation of the \MOR approach for the model is described in detail and results compared to the high-fidelity model are presented. Subsequently, a conclusion is drawn in section~\ref{sec:conclusion}.    

\section{Methodology}
\label{sec:method}
In the following, we present a methodology that can approximate a high-fidelity model for \emph{time-dependent} simulation parameter functions and a fixed initial value without any knowledge of its underlying dynamics. The example considered in this paper is an elastic multibody model. We want to emphasize that the unknown dynamics can also be described by a finite element (\FE)-discretized partial differential equation (\PDE).

\subsection{Model}
\begin{figure}
	\centering
	\includegraphics[trim={0 0.85cm 0 1.04cm}, clip, width=.6\linewidth]{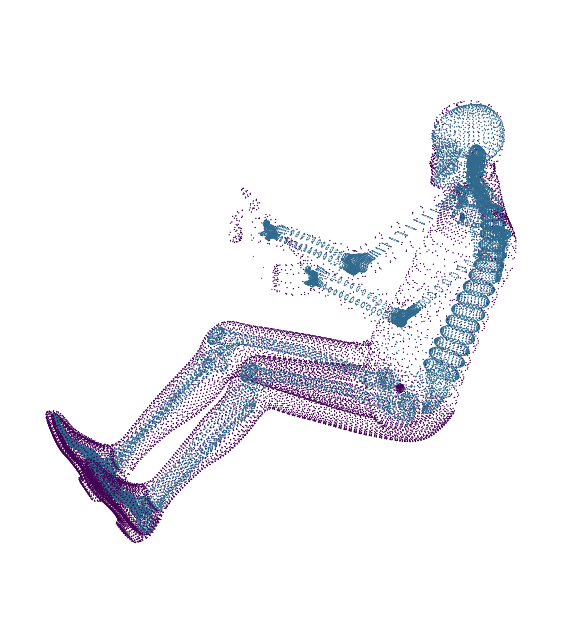}
	\caption{Nodes of the human occupant model.}
	\label{fig:model2}
\end{figure}

The Madymo simulation model consists of a vehicle interior and a reactive human body model which is positioned on the driver seat with its hands on the steering wheel and fastened seat belt, see \cite{TASS19}. For a given driving scenario the vehicle accelerations are applied to the simulation model. The movement of the occupant is measured by recording the nodal displacement during the simulation.
By neglecting all static vehicle interior objects only the relevant nodes that describe the displacement of the human body model remain, see Fig.~\ref{fig:model2}. Further information about model quantities follow in section~\ref{sec:implementation}.
However, the internal calculation rules remain inaccessible.

\subsection{Problem Definition}
Accordingly, only the the high-fidelity \emph{black-box} model 
\begin{align}
\state(t)=\bm{F}(t,\param,\state_{1})
\label{eq:flow map}
\end{align}
is available as basis for the creation of a surrogate model. In this context $t\in\timeInt=\{t_1,...,t_{\nTime}\}$ describes the considered set of discrete time points, $\state\colon\timeInt \to \stateSpace \subseteq \R^{\stateDim}$ represents the state vector describing the system and $\param\colon\timeInt \to \paramSpace \subseteq \R^{\paramDim}$ the simulation parameter function. The solution operator~\eqref{eq:flow map} is known as \emph{flow} map and maps the initial condition~ $\state_1\in\R^{\stateDim}$ and the past values of the time-dependent simulation parameter trajectory~$\param$ to the solution at time $t\geq t_1$. %, \ t_{i+1}=t_i+\Delta t 

From a mathematical point of view, the task of approximating the flow map~\eqref{eq:flow map} can be formulated as optimization problem
\begin{align}
	\min_{\transition} \quad & \underset{\substack{\\\param\in\paramFuncSpace}}{\text{mean}} \quad	
	J(\state(\cdot),\stateApprox(\cdot))	&  
	\label{eq:problem definition}\\
	\text{s.t.} \quad	&\stateApprox(t_{i+1}) = \transition(\stateApprox(t_{i}), \param(t_i)), &i\in \indSpace \setminus \{\nTime \}  \nonumber \\
						&{\state}(t) = \Flow(t,\param,\state_1), 	&t\in\timeInt \nonumber\\
						&\stateApprox_1 = \stateApprox(t_1) = \state(t_1) = \state_1, \nonumber &
\end{align}
where $\indSpace:=\{1,...,\nTime\}$ represent the time indices, and $\paramFuncSpace:= \{ \param:\timeInt \to \paramSpace \}$ the admissible set of parameter functions. 
The main purpose of the optimization problem~\eqref{eq:problem definition} is to find a surrogate model~$\transition$ with output $\stateApprox$, which approximates the high-fidelity solution best with respect to the objective function~$J(\state(\cdot),\stateApprox(\cdot))$. Moreover, the mean operator ensures that the optimum is reached for a wide variety of considered parameter trajectories. Note that another goal is to create a surrogate model, which is cheaper to evaluate than the full model. 
To simplify the presentation, it is assumed that all all discrete-time quantities are evaluated at the same time points, i.e on the set~$\timeInt$. Nevertheless, the following methodology can be applied to series of varying lengths as well.
%%%
\subsection{Model Reduction}
\label{sec:model reduction}
The dynamic solution of high-dimensional \FE or elastic multibody systems are often approximately embedded in a low-dimensional space~${\redStateSpace \in\R^\redStateDim}$ with~$r\ll N$. Hence, an approximation in the reduced order space can be beneficial regarding computational complexity and interpretability. One approach to find a suitable low-dimensional description are reduced-basis methods, see \cite{QuarteroniManzoniNegri16}. They aim to capture the original system behavior with a linear combination~$\state(t)\approx\sum_{l=1}^{r}\bar{z}^{(l)}(t) \ \redBasis_l$ of the so-called \emph{reduced basis vectors}~$\{\redBasis_1,...,\redBasis_r\}\subseteq\stateSpace$. The \emph{reduced states} ${\redState(t)=\begin{bmatrix}
	\bar{z}^{(1)}(t) & ... & \bar{z}^{(r)}(t)\end{bmatrix}^T \colon\timeInt \to \redStateSpace \subseteq \R^{\redStateDim}}$ correspond to the reduced basis' coefficients. 

Snapshot-based methods are one class of completely data-driven approaches to obtain the reduced basis functions. 
They rely on a snapshot matrix $\snaps\in\R^{N \times \nSnaps\nTime}$ with
\begin{align*}
	\snaps
	\hspace{-.4mm}=\hspace{-.4mm}
	\begin{bmatrix}
		\state(t_1, \param_1^{(1)}) \
		... \
		\state(t_{\nTime}, \param_{\nTime}^{(1)})\
		...\
		\state(t_1, \param_{1}^{(\nSnaps)})\
		...\
		\state(t_{\nTime}, \param_{\nTime}^{(\nSnaps)})
	\end{bmatrix}
\end{align*} consisting of $\nSnaps$ high-fidelity simulations. \\
Those simulations are based on a discrete and finite parameter set~$\paramSet = \{\param^{(1)},...,\param^{(\nSnaps)}\} \subseteq \paramFuncSpace$ including several relevant load cases. 
It is assumed that the subspace~${\mathcal{S}_{\paramSet}=\text{span}\{\snaps\}}$ approximates the discrete solution manifold~$\mathcal{S}_{\paramFuncSpace}=\{\state(t{,\param}) \ \vert \ t\in\timeInt, \ \param\in\paramFuncSpace\}$ to a satisfying degree as long as enough snapshots are present.

The \POD is one famous snapshot-based method in which a set of ordered orthonormal basis functions $\redBasis_l, \ l=1,...,d$ is determined. These functions can optimally represent the snapshots in the solution manifold with only the first $r$ basis functions, see \cite{GubischVolkwein17}.
Such a set can be computed using the singular value decomposition (\SVD)
\begin{align*}
	\snaps=\bm{U}
	\begin{bmatrix}
		\bm{S} & 0 \\
		0 & 0
	\end{bmatrix}
	\bm{P}^T
	= \bm{U}\bm{\Sigma}\bm{P}^T,
\end{align*}
where the columns of $\bm{U}\in\R^{N \times N}$ are known as left-singular vectors and the columns of $\bm{P}\in\R^{\nTime \times \nTime}$ correspondingly as right-singular vectors. The singular values are stored in decreasing order on the diagonal of~${\bm{S}=\text{diag}(\sigma_1,..., \sigma_d)\in\R^{{d \times d}}}$ which is a part of $\bm{\Sigma}\in\R^{{N \times \nTime}}$. 
Consequently, the column space of $\snaps$ can be represented with $d$ linear independent columns of $\bm{U}$ considering $\bm{\Sigma}\bm{P}^T$ as coefficients.
The optimal $r$-dimensional approximation 
\begin{align}
	{\state}(t)\approx\tilde{{\state}}(t):=\bm{V}\bar{{\state}}(t)=\bm{U}_r\bar{{\state}}(t)
	=\begin{bmatrix} 
		\redBasis_1 & ... & \redBasis_r
	\end{bmatrix}
	\bar{{\state}}(t),
	\label{eq:approximation}
\end{align} 
of the system states with respect to the Frobenius norm is found by truncating $\bm{U}$ yielding $\bm{V}\in\R^{N \times r}$.
%Hence, it is possible to describe the system dynamics in their reduced representation $\redState(t)$ as long as the reconstruction of the system states $\reconstruction\reduction\state(t)$ satisfies the approximation requirements.
%
\subsection{Surrogate Model Architecture}
\label{sec:model architecture}
Once a reduced description of the system states is found, the question arises of how to approximate their trajectories. One possibility is an \emph{iterative} approach such that the respective parameter dependency can be considered at each time step. 
Therefore, the discrete state transition is described by the difference equation $\state(t_{i+1}) = \state(t_i) + \Delta \state(t_i)$ with the state difference $\Delta \state(t_i)\in \stateSpace$. We also transfer this assumption to the reduced dynamics, which allows us to approximate the difference equation in the reduced space and subsequently project it back into the physical one. 
Consequently, the surrogate model should approximate a time stepping scheme and takes the form
\begin{align} 
	\transition(\stateApprox(t_i), \param(t_i)) 
	&=\reconstruction(\reduction\stateApprox(t_i) + \difference(\reduction\stateApprox(t_i), \param(t_i))) \nonumber\\
	&=\reconstruction(\redStateApprox(t_i) + \difference(\redStateApprox(t_i), \param(t_i)))  
	\label{eq:transition} \\
	&=\reconstruction(\redStateApprox(t_i) + \deltaRedStateApprox(t_i)), \nonumber
\end{align}
with an unknown function $\difference\colon\redStateSpace \times \paramSpace \to \Delta\redStateSpace$ approximating 
\begin{align}
	\difference(\redState(t_i), \param(t_i)) = \deltaRedStateApprox(t_i) \approx \deltaRedState(t_i)=\reduction\Delta\state(t_i).
	\label{eq:difference}
\end{align}
The field of \emph{machine learning} (\ML) provides a library of regression algorithms to identify \eqref{eq:difference} solely based on data. In particular, a parametric regression algorithm to find $\tilde{\difference}(\redState(t_i), \param(t_i),\weights)$ is chosen for the identification of the unknown function. In this context, $\weights$ represent the additional parameters of the algorithm. They are optimized during training such that the predicted state difference $\deltaRedStateApprox=\tilde{\difference}(\redState(t_i), \param(t_i),\weights)$ approximates $\deltaRedState$ best with respect to a certain loss function $\mathcal{L}(\deltaRedStateApprox, \deltaRedState)$. A common loss function is the squared error loss 
\begin{align}
	\loss_\text{se}(\deltaRedState, \deltaRedStateApprox) := \frac{1}{r}\sum_{l=1}^{r}(\Delta{\bar{z}}^{(l)} - \Delta\tilde{\bar{z}}^{(l)})^2.
	\label{eq:loss}
\end{align} 
\subsection{Regression via Long Short-Term Memory} 
\label{sec:regression lstm}
The class of recurrent neural networks (\RNNs), to which \emph{long short-term memory networks} (\LSTMs) belong, store information of previous temporal steps in a hidden state using feedback connections, see \cite{Aggarwal2018}. Hence, they are frequently used for the prediction of sequential data and are well suited to solve the time-dependent problem of approximating~$\difference(\redState(t_i), \param(t_i))$. 
In a classical recurrent network, the hidden state~$\NNhiddenState_i\in\R^{\NNhiddenDim}$ and the output~$\NNoutput_i\in\R^{\NNoutputDim}$ at the $i$-th time step are calculated by
\begin{align}
	\NNhiddenState_i=\sigmoid_h(\Weights_{hx}\NNinput_i+\Weights_{hh}\NNhiddenState_{i-1}+\bias_{\NNhiddenState_i})&=\bm{f}^{\NNhiddenState}(\NNinput_i, \NNhiddenState_{i-1}) 	\label{eq:recurrent NN}\\
	\NNoutput_i = \sigmoid_y(\Weights_{yh}\NNhiddenState_i+\bias_{\NNoutput})&=\bm{f}^{\NNoutput}(\NNinput_i, \NNhiddenState_{i-1}),
\end{align}
where~$\NNinput_i\in\R^{\NNinputDim}$ is the input, $\Weights_{hx}\in\R^{\NNhiddenDim\times\NNinputDim}$ represents the input to hidden weight matrix,~$\Weights_{hh}\in\R^{\NNhiddenDim\times\NNhiddenDim}$ the hidden to hidden and~$\Weights_{yh}\in\R^{\NNoutputDim\times\NNhiddenDim}$ the hidden to output weight matrix with biases~$\bias_{\NNhiddenState}\in\R^{\NNhiddenDim}$ and $\bias_{\NNoutput}\in\R^{\NNoutputDim}$. Furthermore, $\sigmoid_h$ and $\sigmoid_y$ are the activation functions of the hidden and output state. 
Each prediction shares the same weights and consequently the same underlying functions are evaluated for all time steps. Moreover, recurrent neural networks can deal with sequences of different lengths due to the recursive form of \eqref{eq:recurrent NN}.

However, classical \RNNs suffer from vanishing or exploding gradients. To mitigate this problem \cite{HochreiterSchmidhuber1997} introduced long short-term memory networks, which use gates to ensure constant error flow through internal states. 
\begin{figure}
	\includegraphics[trim={0 0cm 0 0cm}, clip, width=\linewidth]{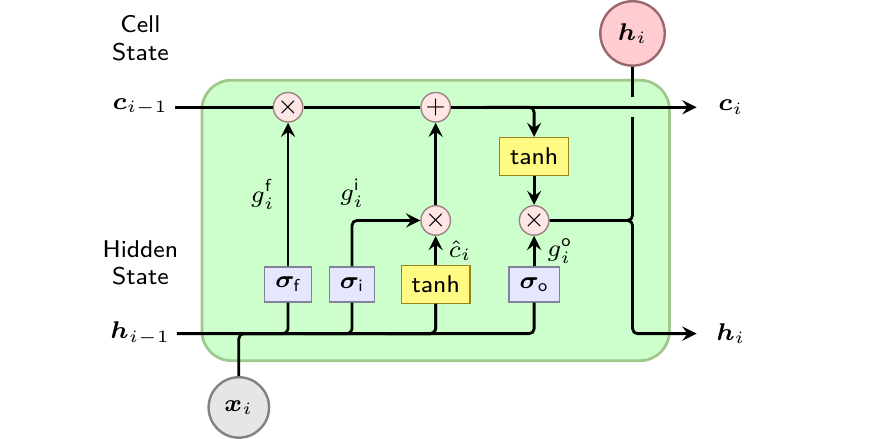}
	\caption{Structure of a \LSTM block with pointwise operations~\protect\tikz \protect\node [circle,draw={rgb,255: red,155; green,130;blue,130},inner sep=-0.5pt,minimum height =.3cm,fill={rgb,255: red,255; green,230; blue,230}] (mux1) at (0,0) {$\cdot$}; 
		and neural network layers for the gates~\protect\tikz 
		\protect\node [rectangle,draw={rgb,255: red,130; green,130; blue,155},inner sep=2pt,fill={rgb,255: red,230; green,230; blue,255}, minimum width=.6cm] (mux1) at (0,0) {$\cdot$};
		with sigmoid activation functions restricting their output to $(0,1)$, i.e. closed and open, and layers with a tanh activation function~\protect\tikz 
		\protect\node [rectangle,draw={rgb,255: red,130; green,120; blue,30},inner sep=2pt,fill={rgb,255: red,255; green,250; blue,130}, minimum width=.6cm] (mux1) at (0,0) {$\cdot$};.
		The gate units $\bm{g}_i^\text{f}$, $\bm{g}_i^\text{i}$ and $\bm{g}_i^\text{o}$ control the data flow. Time-dependencies are propagated through the hidden state.}
	\label{fig:lstm module}
\end{figure}
The structure of a \LSTM block is shown in Fig.~\ref{fig:lstm module} and can be described by 
\begin{align*}
	\bm{g}_i^\text{f}&=\sigmoid_\text{f}(\Weights_\text{f}[\NNhiddenState_{i-1},\NNinput_i]+\bias_\text{f}) 	\hspace{3.7mm}\hat{\bm{c}}_i=\tanh(\Weights_c[\NNhiddenState_{i-1},\NNinput_i]+\bias_c)\\
	\bm{g}_i^\text{i}&=\sigmoid_\text{i}(\Weights_\text{i}[\NNhiddenState_{i-1},\NNinput_i]+\bias_\text{i}) 
	\hspace{4.2mm}\bm{c}_i=\bm{g}_i^\text{f} \bm{c}_{i-1}+\bm{g}_i^\text{i}\hat{\bm{c}}_i \\
	\bm{g}_i^\text{o}&=\sigmoid_h(\Weights_h[\NNhiddenState_{i-1},\NNinput_i]+\bias_h) 
	\hspace{1.5mm}\NNhiddenState_i=\bm{g}_i^\text{o}\tanh(\bm{c}_i)
\end{align*}
with the input $\NNinput_i \in \R^{\NNinputDim}$, the cell state $\bm{c}_i \in \R^{\NNhiddenDim}$ which is used for the long term memory capability, the hidden state $\NNhiddenState_{t-1} \in \R^{n_{h}}$ serving as output of the cell and the weights $\Weights_\text{f}, \Weights_\text{i}, \Weights_c, \Weights_h \in \R^{\NNhiddenDim\times({n_\text{h}+\NNinputDim)}}$ and biases $\bias_\text{f}, \bias_\text{i}, \bias_c, \bias_h \in \R^{n_\text{h}}$.
In this context, the forget gate~$\bm{g}_i^\text{f} \in \R^{\NNhiddenDim\times({n_\text{h}+\NNinputDim)}}$ decides whether any information remains in a cell, the input gate~$\bm{g}_i^\text{i} \in \R^{\NNhiddenDim\times({n_\text{h}+\NNinputDim)}}$ determines to what extent new values flow into a cell and the output gate~$\bm{g}_i^\text{o} \in \R^{\NNhiddenDim\times({n_\text{h}+\NNinputDim)}}$ specifies the information that leaves the cell. All gates use sigmoid activation functions ($\sigmoid_\text{f}$, $\sigmoid_\text{i}$ and $\sigmoid_{h}$).
For further information about \LSTMs refer to \cite{HochreiterSchmidhuber1997}.
Let $\NNhiddenState_i = \bm{f}^h(\NNinput_i,\NNhiddenState_{i-1})$ be the function that maps a current input and previous hidden state to the following one. Then
\begin{align*}
	\NNhiddenState_i &= \bm{f}^h(\NNinput_i,\NNhiddenState_{i-1}) \\ \vspace{-15mm}
	&=\bm{f}^h(\NNinput_i, 
	\underset{\NNhiddenState_{i-1}}{
		\underbrace{
			\bm{f}^h(\NNinput_{i-1}, 
			\bm{f}^h(..., 
			\bm{f}^h(\NNinput_{i-\timeHorizon+1}, 
			\NNhiddenState_{i-\timeHorizon}
			)
			)
			)
	}}
	)\\
	&=:\mathcal{F}(\NNinput_i,\NNinput_{i-1},...,
	\NNinput_{i-\timeHorizon+1}, \NNhiddenState_{i-\timeHorizon})
\end{align*}
calculates the current hidden state based on the $\timeHorizon$ last states by sequential evaluation of $\bm{f}^h$. Please note that the same weights are used in each evaluation. In the last layer of the network the hidden states correspond to its output, i.e $\NNoutput_{i}=\NNhiddenState_{i}=\mathcal{F}(\NNinput_i,\NNinput_{i-1},...,
\NNinput_{i-\timeHorizon+1}, \NNhiddenState_{i-w})$. Hence, the multi time step dependency also transfers to the output and the size of the last layer must be compatible with the output dimension $\R^{\NNoutputDim}$.

In the problem considered in this paper, the state difference $\deltaRedState_{i}(\param^{(j)})$ at the $i$-th time step from the $j$-th simulation serves as network output data $\NNoutput_{ij}$. The input $\NNinput_{ij}$, on the other hand, consists of the previous $\timeHorizon$ reduced basis coefficients $\redState$ along with the corresponding simulation parameters $\param$. This all leads to a dataset
% data set
\begin{align}
	\dataSet = 
	\Bigg(\underset{\NNinput_{ij}\in\R^{(r+n_p)\times \timeHorizon}}
		{\underbrace{ 
		\begin{bmatrix}
			\bar{\state}_{i-\timeHorizon}(\param_{i-\timeHorizon}^{(j)}) 
				\hdots 
			\bar{\state}_{i-1}(\param_{i-1}^{(j)})\\
			\param_{i-\timeHorizon}^{(j)}
				\hdots 
			\param_{i-1}^{(j)}
		\end{bmatrix}}}, \quad
	\underset{\NNoutput_{ij}\in\R^{r}}
		{\underbrace{
		\begin{bmatrix}
			\deltaRedState_{i}(\param^{(j)})
		\end{bmatrix}}}
	\Bigg)^{j\in\mathcal{K}}_{i\in\timeInt}
	\nonumber
%	\label{eq:dataset}
\end{align}
which can be assembled out of the high-fidelity simulation results with $\mathcal{K}=\{1,...,\nSnaps\}$. 
We emphasize that the first $w-1$ data samples per simulation lack predecessors. Based on this dataset, the optimization process of the weights $\weights$ of the \LSTM network \eqref{eq:difference} can be conducted.
%%%%%%%%%%%%%%%%%%%%%%%%%%%%%%%%%%%%%%%%%%%%%%%%%%%%%%%%%%%%%%%%%%
%% 					Problem Reformulation						%%
%%%%%%%%%%%%%%%%%%%%%%%%%%%%%%%%%%%%%%%%%%%%%%%%%%%%%%%%%%%%%%%%%%
\subsection{Problem Reformulation} 
The modeling ideas mentioned so far lead to a reformulation of the initially stated problem definition. By introducing the ideas of \emph{reduction} from section~\ref{sec:model reduction}, transition via state differences from section~\ref{sec:model architecture} and \emph{regression} from section~\ref{sec:regression lstm} along with the corresponding loss function \eqref{eq:loss} into~\eqref{eq:problem definition}, the new problem formulation
\begin{align}
	\min_{\weights} \quad & \underset{\substack{\\\param\in\paramFuncSpace_\nSnaps}}{\text{mean}} \quad	
	\loss_\text{se}(\deltaRedState, \deltaRedStateApprox)	& & 
	\label{eq:problem reformulation}
	\\
	\text{s.t.} \quad &\redState(t_{i+1}) = \redState(t_i) + {\tilde{\difference}(\redState(t_i), \param(t_i), \weights)}, &i&\in\indSpace\setminus\{\nTime\}	\nonumber\\
	&\stateApprox(t) = \redMatrix\redState(t), 				&t&\in\timeInt&  \nonumber\\
	&{\state}(t) = \bm{F}(t,\param,\state_1), 				&t&\in\timeInt&\nonumber\\
	& \redState_1 = \redState(t_1) = \redMatrix^T\state_1. \nonumber
\end{align} 
arises. For practical reasons, only a finite set of parameter functions~$\paramFuncSpace_\nSnaps$ is considered in~\eqref{eq:problem reformulation}. Furthermore, 
it is valid to replace~${J}(\state, \stateApprox)$ in favor of~$\loss_\text{se}(\deltaRedState, \deltaRedStateApprox)$ given that the overall objective function is minimized as the predicted state difference approaches the reference data.
%%%%%%%%%%%%%%%%%%%%%%%%%%%%%%%%%%%%%%%%%%%%%%%%%%%%%%%%%%%%%%%%%%
%% 							ALGORITHM  							%%
%%%%%%%%%%%%%%%%%%%%%%%%%%%%%%%%%%%%%%%%%%%%%%%%%%%%%%%%%%%%%%%%%%
\subsection{Algorithm}

All mentioned ideas are combined into one algorithm which will be referred to as \POD-\LSTM in the following. 
It can be divided into two parts: (i) an offline process to solve~\eqref{eq:problem reformulation}, summarized in Algorithm~\ref{alg:Offline}, in which the high-fidelity model is evaluated, the reduction matrix $\redMatrix$ is calculated, the dataset $\dataSet$ is assembled and the training of the \LSTM network is conducted and (ii) an online process, in which the dynamics are approximated for new parameter trajectories summarized in Algorithm~\ref{alg:Online}.
 
During the online part, the difference $\Delta\redStateApprox$ is predicted based on the reduced initial state $\redState_1$ and the corresponding parameter $\param_1$. Hereafter, the resulting reduced state $\redStateApprox$ is calculated, which in turn serves as input for the next prediction. Hence, the first predictions lack predecessors. Thus, the ability of \LSTMs to deal with inputs of variable length is of high importance.
\begin{algorithm}
	\caption{\textsf{POD-LSTM Offline }}
	\textbf{Input}: parameter domain $\paramFuncSpace$, amount of simulations $\nSnaps$\\
	\textbf{Output}: reduction matrix $\redMatrix$, regression model $\tilde{\difference}$
	\begin{algorithmic}[1]
		\State create parameter set $\paramFuncSpace_\nSnaps=\{\param^{(1)}(t),...,\param^{(\nSnaps)}(t)\}$
		\For {$j=1,...,\nSnaps$}
			\State evaluate black-box model $\Flow(t,\param^{(j)}(t),\state_1)$
		\EndFor
		\State assemble snapshot matrix $\snaps$
		\State calculate reduction matrix $\redMatrix$ via \POD 
		\State assemble dataset $\dataSet$
		\State train \LSTM network $\tilde{\difference}(\redState(t_i), \param(t_i), \weights)$
	\end{algorithmic}
	\label{alg:Offline}
\end{algorithm}
\begin{algorithm}
	\caption{\textsf{\POD-\LSTM Online }}
	\textbf{Input}: initial state $\state_1\in\stateSpace$, amount of timesteps $\nTime$, \\ \phantom{asd} parameter $\param(t)\in\paramFuncSpace$\\
	\textbf{Output}: approximated state trajectory $\stateApprox(t)$ \nonumber
	\begin{algorithmic}[1]
		\State reduce initial state $\redState[1] = \redMatrix^T\state_1$
		\For {$t=1,...,\nTime-1$}
			\If{$t-w<1$}
			\State $i = 1$
			\Else
			\State $i = t-w$
			\EndIf
			\State $\redStateApprox[t+1] = \redStateApprox[t] + \tilde{\difference}(\redStateApprox[i:t],\param[i:t])$
		\EndFor
		\State project into full space $\stateApprox[:]=\redMatrix\redStateApprox[:]$
		\State \textbf{return} $\stateApprox$
	\end{algorithmic}
	\label{alg:Online}
\end{algorithm}
\section{Implementation \& Results}
\label{sec:implementation}
As first step of Algorithm~\ref{alg:Offline}, a set of $107$ parameter trajectories $\paramFuncSpace_{107}=\{\param^{(1)}(t),...,\param^{(107)}(t)\}$ is generated for the simulations with the human occupant model. Each trajectory $\param^{(i)}(t)\in\R^3$ includes accelerations in all three Cartesian coordinate directions. Subsequently, the high-fidelity occupant model is evaluated for each trajectory. 
The generated simulation results are randomly split into training, validation and test data. The training data consists out of $\nSnaps=90$ full simulations resulting in $n_\text{sample}=10478$ samples, whereas the validation set includes 11 and the test data 6 simulations.
All simulation results are sampled every $0.025$\,s and contain the displacements $\nodeDisp_n(t)=\begin{bmatrix}\nodeDisp_{1}(t) & \nodeDisp_{2}(t) & \nodeDisp_{3}(t)\end{bmatrix}\in\R^{n_\text{node}\times 3}$ of all $n_\text{node}=\nodeDim$ nodes with respect to a reference configuration. In this regard, $\nodeDisp_{1}(t)$ represent the displacement in the first coordinate direction and $\nodeDisp_{2}(t)$ and $\nodeDisp_{3}(t)$ the displacement in the corresponding second and third one. 
For further processing, the displacements are vectorized yielding the system states 
$\state(t)=\text{vec}(\nodeDisp(t))\in\R^{\stateDim}$ with $\stateDim=\dofs$.

\begin{figure}
	\includegraphics[width=.5\textwidth]{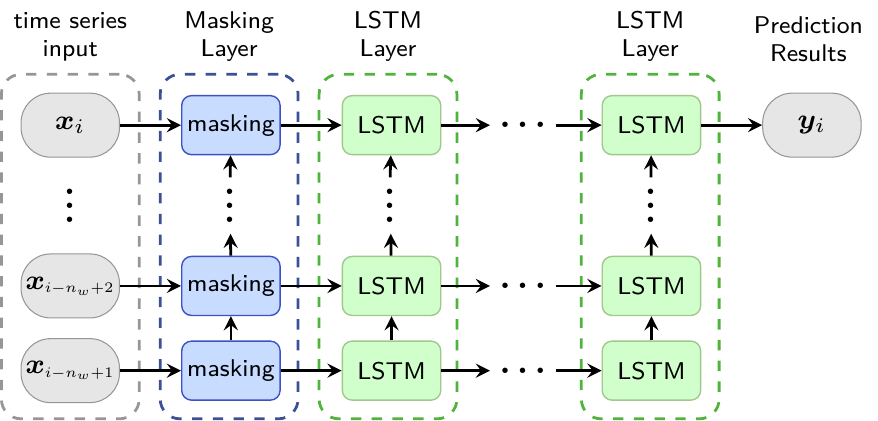}
	\caption{Network architecture including a masking layer to deal with inputs of varying length, two subsequent \LSTM layers to capture time dependencies and a closing fully connected layer. For reasons of clarity, the first three layers are unfolded in time.}
	\label{fig:Network Architecture}
\end{figure}

Based on the training data, the snapshot matrix $\snaps\in\R^{\stateDim \times n_\text{sample}}$ is assembled and the remaining \POD-\LSTM  offline algorithm is evaluated.
The reduced system size is chosen to be $\nRed=30$.
Tensorflow\footnote{https://www.tensorflow.org/} is used for the implementation of the neural network. The architecture consists of a masking layer followed by four \LSTM layers as shown in Fig.~\ref{fig:Network Architecture}. The masking layer ensures comparability with samples of different length. The first \LSTM layers consist of $256$ units each, while the last one only contains $r=30$ units to ensure compatibility with the output dimension. It applies for all \LSTM layers, that the previous $8$ time steps are considered for each prediction. The network is trained over $150$ epochs using the RMSprop optimizer with a learning rate of $0.001$ and a batch size of 5. 
The network weights $\weights$ that resulted in the lowest loss for the validation set during training are selected for the following evaluation. All mentioned hyperparameter are fitted using a grid-search on the validation set. 

\subsection{Error Quantities}
The test data consists of simulations that are not used during the offline stage, neither for the calculation of the reduced basis using \POD nor during the training of the network. For each simulation, the \POD-\LSTM online algorithm is evaluated to receive approximated state trajectories. For the comparison against the reference solution obtained by the high-fidelity model, time-dependent relative error scores 
\begin{align}
	s({\bm{\zeta}}, \tilde{{\bm{\zeta}}}, t)=1-\frac{\Vert {\bm{\zeta}}(t) - \tilde{{\bm{\zeta}}}(t) \Vert_2}{\Vert \mathbb{\bm{\zeta}}(t) \Vert_2}
	\label{eq:error}
\end{align}
are introduced.
A value of $s=1$ corresponds to a perfect approximation, whereas the quantity can become negative since the surrogate model can be arbitrarily bad.
Furthermore, the mean of the relative score per simulation $\hat{s}({\bm{\zeta}}, \tilde{{\bm{\zeta}}})=\frac{1}{\nTime}\sum_{i=1}^{\nTime}s({\bm{\zeta}}, \tilde{{\bm{\zeta}}}, t_i)$ serves as overall metric.

Substituting a test state 
$\state_\text{test}(t)=\state(t, \param_{\text{test}})$
and its reconstructed equivalent obtained via \POD into \eqref{eq:error} yields the
\emph{reconstruction score} $\eRec(\state_\text{test}, \redMatrix\redMatrix^T\state_\text{test}, t)=s(\state_\text{test}, \redMatrix\redMatrix^T\state_\text{test}, t)$. On the contrary, the \emph{regression score} $\eRegr(\redState_\text{test}, \redStateApprox_\text{test}, t)=s(\redState_\text{test}, \redStateApprox_\text{test}, t)$ compares the reduced basis coefficients $\redState_\text{test}(t)$  with their approximation obtained from the \LSTM network. This score does not only result from the networks' prediction error but also from error propagation arising from the iterative nature of the approach. The next considered quantity $\eApprox(\state_\text{test}, \stateApprox_\text{test}, t)=s(\state_\text{test}, \stateApprox_\text{test}, t)$ represents the \emph{approximation score} and thus the overall approximation quality achieved by the \POD-\LSTM algorithm.

In order to introduce a more tangible metric, the Euclidean node distance 
\begin{align*}
\eDist(t, \bm{q}^m, \tilde{\bm{q}}^m) = \lVert \bm{q}^m(t)-\tilde{\bm{q}}^m(t) \rVert_2 
\end{align*} 
between a high-fidelity simulation and its approximation is used. It compares the displacement of the $m$-th node $\bm{q}^m(t)$ with its approximation $\tilde{\bm{q}}^m(t)$. 
Furthermore, the maximum distance over all nodes and time steps 
\begin{align*}
	\eDistMax = \max_{t\in\timeInt}\max_{m\in\{1,...,n_\text{node}\}} \eDist(t, \bm{q}^m, \tilde{\bm{q}}^m)
\end{align*} 
describes the overall error during a simulation.
The last considered criterion measures the model's computation time by calculating the ratio between simulation time $t_\text{cpu}$ and simulated time, i.e. 
\begin{align*}
	\tRatio=\frac{t_\text{cpu}}{t_\nTime-t_1}.
\end{align*} 
A value of $\tRatio\leq1$ thus is equivalent to real-time capability. It must be mentioned that only the online computation time is included in this metric and the offline computation time is neglected.

\subsection{Simulation Results}
In this section, the performance of the surrogate model is validated for 6 test simulations. All high-fidelity and surrogate simulations, are conducted on an Intel Core i7-6700 CPU with 32\,GB of RAM. The simulations vary in length and consist out of 74 up to 164 test samples resulting in a total amount of 600 test samples. 
The results are summarized in Table~\ref{tab:perf}. 
\begin{table}
	\setlength{\tabcolsep}{5pt}
	\centering
	\caption{Error quantities of 6 test simulations.}
	\label{tab:perf}
	\begin{tabular}[c]{c c c c c c c c}
		\toprule
		&sim 1		& sim 2 	& sim 3 	& sim 4	& sim 5	& sim 6	& mean\\ 
		\midrule
		$\hat{s}_\text{regr}$\,(-) 		&0.948 & 0.963 & 0.856 & 0.839 & 0.961 & 0.955 & 0.920\\ 
		$\hat{s}_\text{appr}$\,(-) 		& 0.947 & 0.961	& 0.855 & 0.838	 & 0.960 & 0.953 & 0.919\\ 
		$\hat{s}_\text{appr}^\text{1}$\,(-) 		 & 0.969	& 0.977 & 0.909	 & 0.922 & 0.977 & 0.981 & 0.956\\ 
		$\hat{s}_\text{rec}$\,(-)		& 0.991 & 0.993 & 0.990 & 0.991 & 0.992 & 0.992 & 0.992\\ 
		$\eDistMax$\,(cm)	& 0.285 & 0.407 & 0.642 & 0.659 & 0.382 & 0.989 & 0.561\\ 
		$\tRatio$\,(-) & 0.982 & 0.989 & 0.984 & 0.975 & 0.990 & 0.990 & 0.985\\ 
		\bottomrule
	\end{tabular} 
\end{table}

It is noticeable that the reconstruction score is almost optimal for all simulations, what confirms the assumption to approximate the dynamics in reduced space. Furthermore, the regression $\eRegr$ and approximation error $\eApprox$ are almost identical. From this it can be concluded that the error within the iterative regression is the decisive factor for the overall approximation quality.

At first glance, there seems to be large differences in the error metrics between the simulations. However, considering only the first second of each simulatio,n as done in Table~\ref{tab:perf} with $\hat{e}_\text{appr}^\text{1}$, it becomes obvious that the main reason is the simulation time. Due to the iterative approach, deviations are propagated through the entire simulation. Hence, once a deviation occurs it continuously contributes to the further error. Nevertheless, the maximum node distance $\eDist$ among all simulations is still below 1\,cm and the surrogate model outperforms its high-fidelity equivalent regarding computation time. 
While the latter yields a value of $\tRatio=2271.70$ in average, the surrogate reaches real-time capability with an average value of $\tRatio=0.985$, so that the surrogate model is more than 2300 times faster.

\begin{figure}
	\includegraphics[trim={0 0.1cm 0 0cm}, clip, width=\linewidth]{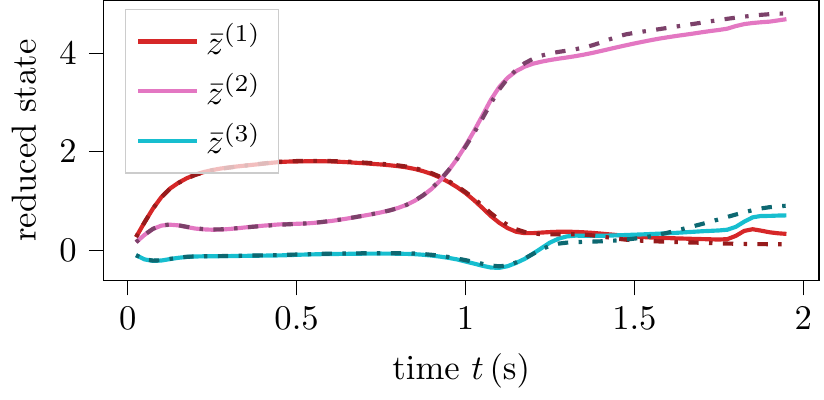}
	\caption{Trajectories of the first three reduced basis' coefficients of one test simulation. The reference coordinates are drawn as solid line~\protect\tikz[baseline=-.75ex]\protect\draw [thick] (0,0) -- (0.4,0);, whereas their approximation is drawn as dash dotted line~\protect\tikz[baseline=-.75ex]\protect\draw [thick,dash dot] (0,0) -- (0.42,0);.}
	\label{fig:regression error}
\end{figure}

For a more detailed analysis, we present the results of simulation 6, which reaches an intermediate approximation quality but exhibits the highest input dynamics.
Therefore, the trajectories of the reduced basis' coefficients $\redState$ are compared with their approximations $\redStateApprox$ in Fig.~\ref{fig:regression error}.
For the sake of clarity, we decided to present only the coefficients of the first three reduced modes since these have the highest relevance. As displayed, the predictions of the \POD-\LSTM algorithm follow the reference trajectories quite well for the longest time. However, it can be observed again that once a deviation occurs, it persists.
%  until the end of the simulation due to the itaretive nature

\begin{figure}
	\includegraphics[trim={0 0.1cm 0 0.1cm}, clip, width=\linewidth]{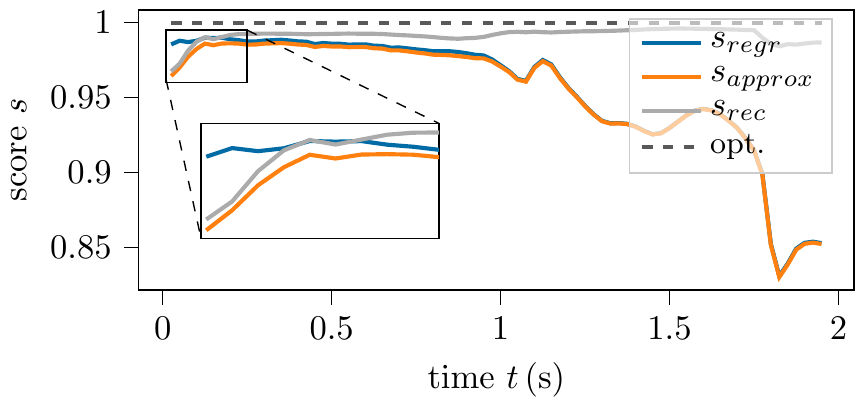}
	\caption{Relative error scores of one test simulation.}
	\label{fig:rel errors}
\end{figure}

Further aspects are highlighted in Fig.~\ref{fig:rel errors}, where the reconstruction score~$\eRec$, approximation score~$\eApprox$ and regression score~$\eRegr$ are shown.
It is observable that the approximation score suffers from both regression and reconstruction. During the first $0.2$\,s, the reconstruction error is large by its standards and thus contributes to the approximation error above average. However, for the rest of the simulation it is close to the optimum and influences the approximation score only marginally.
Consequently, the course of $\eApprox$ corresponds to that one of $\eRegr$ to a large extent allowing conclusions to be drawn directly from the approximation in the reduced space. 

\begin{figure}
	\includegraphics[trim={0 0.1cm 0 0.2cm}, clip, width=\linewidth]{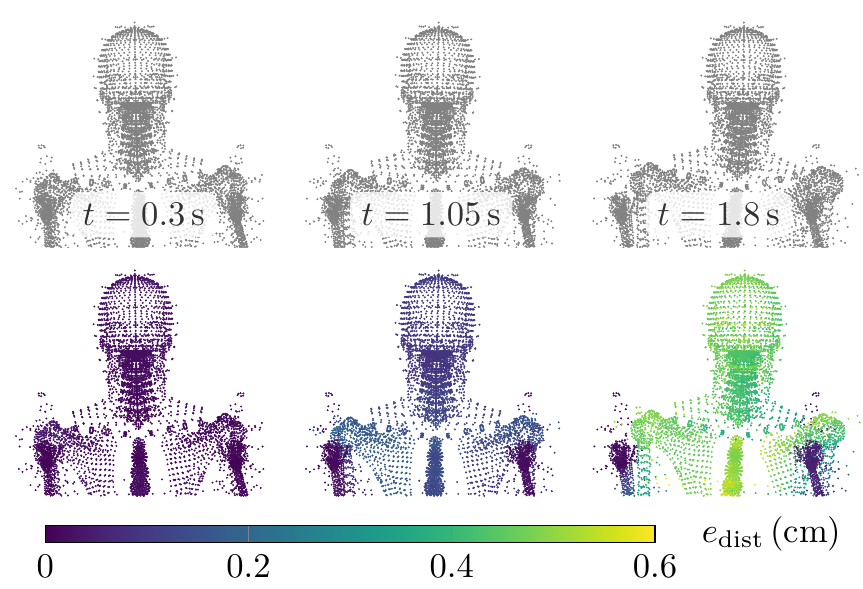}
	\caption{Motion of the occupant's upper body at certain time points. The reference solution is shown in gray in the upper half of the figure, whereas the nodes of the approximation in the bottom half are colored with respect to their individual node distance $\eDist$.}
	\label{fig:eucl dist}
\end{figure}

As final comparison, a visualization of the occupant's motion and its approximation for certain time steps is recorded in Fig.~\ref{fig:eucl dist}. 
During the first time steps, almost none differences can be perceived. Only towards the end, a quantitative deviation of the upper body can be observed. The motion sequence is nevertheless well reproduced over the complete course.
Accordingly, the surrogate model approximates the high-fidelity occupant model to a high degree requiring much less computational effort.

\section{Conclusion}
\label{sec:conclusion}
The presented combination of model order reduction and machine learning is able to create a surrogate model for an accelerated human occupant model. For simulations with short duration, even complex time-dependencies could be captured while simultaneously a remarkable speedup of computation time is achieved. 
This allows not only the use in real-time scenarios and scanning of multi-parameter spaces, but also the implementation on low-end hardware. Moreover, the surrogate model can be deployed in further simulations environments without adding much computational effort. 
However, it became apparent that the iterative architecture can lead to error propagation and thus non-negligible deviations to the reference model. Accordingly, it must be a goal of future work to enhance the long term predictions. Either by updating the surrogate model after a certain period of time based on sensor data or by increasing the robustness of the regression algorithm.
%
%\begin{ack}
%Place acknowledgments here.
%\end{ack}
%
%\bibliography{ifacconf}             % bib file to produce the bibliography                                                    % with bibtex (preferred)
%\ifwindows
%	\bibliography{C:/Users/jonas/bwSyncShare/02_Dev/12_itm_literature/ITM_Literatur.bib}%
%\else
%	\bibliography{/scratch/tmp/jkneifl/02_Dev/12_itm_literature/ITM_Literatur.bib}%
%\fi   

%\appendix
%\section{A summary of Latin grammar}    % Each appendix must have a short title.
%\section{Some Latin vocabulary}              % sections and subsections are supported  
%                                                                         % in the appendices.
\end{document}